
\magnification=1100
\overfullrule=0pt
\hsize=15truecm
\vsize=20.5truecm
\baselineskip=21pt
 \baselineskip=12pt



\newfam\msyfam
\font\tenmsy=msbm10
\font\sevenmsy=msbm7
\font\fivemsy=msbm5

\textfont\msyfam=\tenmsy
\scriptfont\msyfam=\sevenmsy
\scriptscriptfont\msyfam=\fivemsy

\def\Bbb{\fam\msyfam\tenmsy}


\def\C{{\Bbb C}}

\def\F{{\Bbb F}}

\def\P{{\Bbb P}}

\def\Z{{\Bbb Z}}


\def\EE{{\cal E}}

\def\HH{{\cal H}}


\def\dim{\hbox{dim}\,}


\def\map{\longrightarrow}

\def\commadots{\ ,\ldots ,\ }

\def\half{{1\over 2}}

\def\section#1
   {\goodbreak\bigskip
   {\noindent\hglue-1.5pt\bf#1}
   \vglue2pt
   \par\nobreak\noindent
}


\font\medbold=cmbx10 scaled 1100

\def\Heading#1{
  \noindent
  {\medbold #1 }
}

\def\heading#1{
  \noindent
  {\bf #1 }
}

\def\address{
\medskip

\begingroup
\obeylines
\parskip=0pt
\parindent=0pt
\baselineskip=12pt
\bf

}

\def\endaddress{\endgroup}

\def\abstract#1//{
\bigskip
\begingroup
\baselineskip=11pt
\noindent{\bf Abstract.}
\noindent{\rm#1}
\endgroup
}

\def\frtitle#1//{
\bigskip
\noindent
{\bf #1}
}

\def\frabstract#1//{
\bigskip
\begingroup
\baselineskip=11pt
\noindent
{\bf R\'esum\'e.}
\noindent
{\rm#1}
\endgroup
}

\def\proclaim#1{
  \goodbreak
  \medskip
  \begingroup
  \noindent
  {\bf #1}
  \sl
  }
\def\endproclaim{
  \endgroup
}

\def\bibliography{
  \bigskip
  \centerline{{\medbold References }}
  \bigskip
  \begingroup
  \parindent1cm
  \parskip4pt
  \baselineskip=11pt
  }

\def\bi#1{%
\item{[#1]}%
}

\def\endbibliography{
  \endgroup
  }

\def\cite#1{[#1]}

\def\ip#1/#2/{\left( \; #1,\; #2\; \right)\;}
\def\hip#1/#2/{h(#1,\, #2 )\;}
\def\emph#1/{{\sl #1}}


\footnote{}{First author partially supported by an NSF postdoctoral fellowship.
 Second
and third authors partially supported by NSF grant DMS 9625463}

\heading{G\'eom\'etrie alg\'ebrique/Algebraic geometry}
\bigskip

\Heading{A Complex Hyperbolic Structure for Moduli of Cubic Surfaces}
\bigskip

\heading{Daniel Allcock, James A. Carlson, and Domingo Toledo}

\address

Department of Mathematics
University of Utah
Salt Lake City, Utah, USA.
E-mail: allcock, carlson and toledo at math.utah.edu

\endaddress

\abstract   We show that the moduli space $M$ of marked cubic surfaces is
biholomorphic to $(B^4 - \HH)/\Gamma_0$ where $B^4$ is complex hyperbolic
four-space,
where
$\Gamma_0$ is a specific group generated by complex reflections,
and where $\HH$ is the union of reflection hyperplanes for $\Gamma_0$. Thus $M$
has
a complex hyperbolic structure, i.e., an (incomplete) metric of constant
holomorphic sectional curvature.
//

\frtitle Une structure hyperbolique complexe pour les modules des surfaces
cubiques//

\frabstract    Nous montrons que l'espace des modules $M$ des
surfaces cubiques marqu\'ees est biholomorphe \`a $(B^4 - \HH)/\Gamma_0$
o\'u $B^4$ est l'espace complexe hyperbolique de dimension quatre,
o\'u $\Gamma_0$ est un groupe sp\'ecifique g\'en\'er\'e par des r\'eflections
complexes, et o\'u $\HH$ est l'union de l'ensemble d'hyperplans de r\'eflection
de $\Gamma_0$. Donc $M$ admet une structure hyperbolique complexe, c'est \`a
dire
une m\'etrique (incompl\`ete) de courbure holomorphe sectionnelle constante.
//

\bigskip
\heading{Version fran\c caise abr\'eg\'ee}

A une surface cubique (marqu\'ee) correspond une vari\'et\'e cubique de
dimension trois
(marqu\'ee), \`a savoir le rev\^etement de $\P^3$ ramifi\'e le long de la
surface.
L'application des p\'eriodes $f$ pour ces vari\'et\'es de
dimension trois est d\'efinie sur l'espace des modules $M$ des cubiques
marqu\'ees, et
cette application $f$ prend ses valeurs dans une quotient de la boule unitaire
dans $\C^4$
par l'action du groupe de monodromie projective.  Ce groupe $\Gamma_0$ est
g\'en\'er\'e par
des r\'eflections complexes dans un ensemble d'hyperplans dont nous notons la
r\'eunion par
$\HH$.  Alors nous avons le resultat suivant:

\proclaim{Th\'eor\`eme.} L'application des p\'eriodes d\'efinit une
biholomorphisme
$$
  f: M \map \left(B^4 - \HH \right)/\Gamma_0.
$$
\endproclaim

\noindent  De ce th\'eor\`eme on obtient des r\'esultats sur la structure
m\'etrique
de $M$ et sur son groupe fondamental:

\proclaim{Corollaires.}  (1) L'espace $M$ des modules des surfaces cubiques
marqu\'ees
admet une structure hyperbolique complexe: une m\'etrique (incompl\`ete) de
courbure
holomorphe sectionnelle constante. (2) Le groupe fondamental de $M$ contient un
sous-groupe
normale qui n'est pas de g\'en\'eration finie. (3)  Le groupe fondamental de
$M$ n'est pas
un r\'eseau dans un groupe de Lie semisimple.
\endproclaim

\medskip
\noindent{\bf Remarques.} (1) Nos m\'ethodes montrent aussi que la
compl\'et\'ee m\'etrique
de $(B^4 - \HH)/\Gamma_0$ est l'orbifolde $B^4/\Gamma_0$,
isomorphe \`a l'espace des modules des surfaces cubiques marqu\'ees stables.
(2)
R\'ecemment E. Looijenga a trouv\'e une pr\'esentation remarkable du groupe
fondamental
orbifolde de l'espace des modules des cubiques lisses non-marqu\'ees.

Afin de pr\'eciser la notion de surface cubique lisse marqu\'ee, fixons un
r\'eseau $L$,
un  $\Z$-module libre avec une base $e_0 \commadots e_6$ qui est muni de la
forme
quadratique telle que la base soit orthogonale et telle que $\ip e_0/e_0/ = 1$,
$\ip e_k/e_k/ = -1$  pour $k > 0$.  Soit $\eta = 3e_0 - (e_1 + \cdots + e_6)$.
Alors
une \emph surface cubique marqu\'ee/ est compos\'ee d'une surface cubique lisse
$S$ et
d'une isom\'etrie $\psi: L \map H^2(S,\Z)$ qui envoie $\eta$ sur la classe d'un
hyperplan.    L'ensemble $M$ de classes d'isomorphisme des surfaces cubiques
marqu\'ees
porte la structure d'une vari\'et\'e et de plus est un \'espace de modules
fines.
Une construction de cette \'espace a \'et\'e donn\'ee dans
\cite{9}, o\'u on trouve aussi une compactification lisse
$C$ de $M$ telle que les points de $C - M$ forment un diviseur \`a croisements
normaux.

Pour d\'efinir le groupe $\Gamma_0$, soit $\EE$ l'anneau des entiers
d'Eisenstein
$\Z[\omega]$ o\'u
$\omega = (-1 + \sqrt{-3})/2$ est une troisi\`eme racine d'unit\'e, et
consid\'erons le
produit Cartesien
$\EE^5$  muni d'une forme hermitienne
$\hip v/w/ = -v_1\bar w_1 + v_2\bar w_2 + v_3\bar w_3 + v_4\bar w_4 + v_5\bar
w_5$.  Alors
$(\EE^5,h)$ est l'unique r\'eseau autodual sur les entiers d'Eisenstein
qui est de signature $(4,1)$. Donc $Aut(\EE^5,h)$ est un sous-groupe discr\`et
du groupe
unitaire $U(h)$, qui agit sur $B^4 =\{\; \ell \in \P^4\; :\;  h|\ell < 0 \;\}$.
 Notons que
$\EE/\sqrt{-3}\EE \cong \F_3$ est un corps de trois el\'ements et notons aussi
qu'il y a
un homomorphisme naturel
$Aut(\EE^5,h) \map Aut(\F_3^5,q)$ o\'u $q$ est la forme quadratique obtenue par
r\'eduction de $h$ modulo $\sqrt{-3}$.  Notons par ``$P$'' la projectivisation,
et
d\'efinissons un groupe $\Gamma_0$ d'automorphismes de $B^4$ par la
suite exacte
$$
  1 \map \Gamma_0 \map PAut(\EE^5,h) \map PAut(\F_3^5,q) \map 1.
$$
Ce groupe est le groupe discr\`et du th\'eor\`eme principal.  Les hyperplans de
$\HH$ sont
d\'efinis par les \'equations $h(x,v) = 0$ pour des vecteurs $v$ dans $\EE^5$
avec $ h(v) =
1$. Notons aussi que $PAut(\F_3^5,q)$ est isomorphe
au groupe de Weyl du r\'eseau
$E_6$.

\section{1.  Main results}

To a (marked) cubic surface corresponds a (marked) cubic threefold defined as
the
triple cover of $\P^3$ ramified along the surface.  The period map $f$ for
these threefolds
is defined on the moduli space $M$ of marked cubic surfaces and takes its
values in the
quotient of the unit ball in $\C^4$ by the action of the projective monodromy
group. This
group $\Gamma_0$ is generated by complex reflections in a set of hyperplanes
whose union we
denote by $\HH$.  Then we have the following result:

\proclaim{Theorem.} The period map defines a biholomorphism
$$f: M \map \left(B^4 - \HH\right)/\Gamma_0.$$
\endproclaim

\noindent  From this identification we obtain
results on the metric structure and the fundamental group:

\proclaim{Corollaries.}  (1) The moduli space of marked cubic surfaces carries
a
complex hyperbolic structure: an (incomplete) metric of constant holomorphic
sectional curvature. (2) The fundamental group of the space of marked cubic
surfaces contains a normal subgroup which is not finitely generated. (3) The
fundamental group of the space of marked cubic surfaces  is not a lattice in a
semisimple Lie group.
\endproclaim

\medskip
\noindent{\bf Remarks.} (1) Our methods also show that the metric completion of
$(B^4 -
\HH)/\Gamma_0$ is the complex hyperbolic orbifold $B^4/\Gamma_0$, which is
isomorphic to
the moduli space of marked stable cubic surfaces.  (2)  Recently E. Looijenga
found a
remarkable presentation of the orbifold fundamental group of the moduli space
of smooth
unmarked cubic surfaces.

To make precise the notion of smooth marked cubic surface, fix the lattice
$L$ to be the free $\Z$-module with basis $e_0 \commadots e_6$
endowed with the quadratic form for which the given
basis is orthogonal and such that $\ip e_0/e_0/ = 1$, $\ip e_k/e_k/ = -1$
for $k > 0$.  Let $\eta = 3e_0 - (e_1 + \cdots + e_6)$.  Then a
\emph marked cubic surface/ consists of a smooth cubic surface $S$ and
an isometry $\psi: L \map H^2(S,\Z)$ which carries $\eta$ to the hyperplane
class.    The set $M$ of isomorphism classes of marked cubic surfaces
has the structure of a variety and is a fine moduli space.  A
construction of it is described in \cite{9}, and a smooth compactification
$C$ is given for which the points of $C - M$ constitute a normal crossing
divisor.

To define the group $\Gamma_0$, let $\EE$ denote the ring of Eisenstein
integers
$\Z[\omega]$ where
$\omega = (-1 + \sqrt{-3})/2$ is a cube root of unity, and consider the
Cartesian product
$\EE^5$  endowed with the hermitian inner product
$\hip v/w/ = -v_1\bar w_1 + v_2\bar w_2 + v_3\bar w_3 + v_4\bar w_4 + v_5\bar
w_5$.  Then
$(\EE^5,h)$ is the unique self-dual lattice over the Eisenstein
integers with signature $(4,1)$. Thus $Aut(\EE^5,h)$ is a discrete subgroup of
the unitary group $U(h)$, which acts on $B^4 =\{\; \ell \in \P^4\; :\;  h|\ell
< 0 \;\}$.
Observe that
$\EE/\sqrt{-3}\EE
\cong
\F_3$ is a field of three elements and that there is a natural homomorphism
$Aut(\EE^5,h) \map Aut(\F_3^5,q)$ where $q$ is the quadratic form obtained by
reduction of $h$ modulo $\sqrt{-3}$.  Let ``$P$'' denote projectivization, and
define a group $\Gamma_0$ of automorphisms of $B^4$ by
$$
  1 \map \Gamma_0 \map PAut(\EE^5,h) \map PAut(\F_3^5,q) \map 1.
$$
This is the discrete group of the main theorem.  The hyperplanes of $\HH$ are
defined by the equations $h(x,v) = 0$ for vectors $v$ in $\EE^5$ with $ h(v) =
1 $. Note
that $PAut(\F_3^5,q)$ is isomorphic to the Weyl group of the
$E_6$ lattice.

\section{2.  Construction of a period mapping}

To construct the period mapping, we examine in detail the Hodge structures
for the cubic threefolds. The underlying lattice $H^3(T,\Z)$ is
ten-dimensional,
carries a unimodular symplectic form $\Omega$, and admits a Hodge decomposition
of the form $H^3(T,\C) = H^{2,1} \oplus H^{1,2}$.  Choose a generator
$\sigma$ for the group of automorphisms of $T$ over $\P^3$, and note that it
operates
without fixed points on
$H^3(T,\Z)$.   This action gives $H^3(T,\Z)$ the structure of a
five-dimensional
module over the Eisenstein integers.  It carries a hermitian form
$$
  \hip x/y/ = \half( \Omega((\sigma - \sigma^{-1})x,\, y)
                      + (\omega - \omega^{-1}) \Omega( x,\, y) )
$$
which is unimodular and of signature $(4,1)$.

Now consider the quotient module $H^3(T,\Z)/(1-\omega)H^3(T,\Z)$ and observe
that it
can be identified isometrically with $(\F_3^5,q)$.   We define a marking of $T$
to be choice of such an isometry, and we claim that a marking of a cubic
surface
determines a marking of the corresponding threefold. Indeed, if $\gamma$ is a
primitive
two-dimensional homology class on $S$  then it is the boundary of a three-chain
$\Gamma$
on $T$. Since $\Gamma$ and $\sigma \Gamma$ have the same boundary, the
three-chain
$c(\gamma) = (1-\sigma)\Gamma$ is a cycle. However, it is well-defined only up
to addition of elements
$(1-\sigma)\Delta$ where $\Delta$ is a three-cycle on $T$.  Thus a homomorphism
$$
  c: H_2^{prim}(S,\Z) \map H_3(T,\Z)/(1-\sigma)
$$
is defined.  Since a marking of $S$ can be viewed as a basis of
$H_2^{prim}(S,\Z)$, application of $c$ to the basis elements defines a basis of
$H_3(T,\Z)/(1-\sigma)$, and this gives the required marking of the threefold.

The action of $\sigma$ decomposes $H^3(T,\C)$ into eigenspaces $H^3_\lambda$
where $\lambda$ varies over the primitive cube roots of unity.  Because
$\sigma$ is
holomorphic, the decomposition is compatible with the Hodge decomposition and
one has
$$
  H^3_\omega = H^{2,1}_\omega \oplus H^{1,2}_\omega
   \qquad
  H^3_{\bar\omega} = H^{2,1}_{\bar\omega} \oplus H^{1,2}_{\bar\omega} .
$$
The dimensions of the Hodge components can be found with the help of Griffiths'
Poincar\'e
residue calculus \cite{5}. Details for this case are found in \cite{3}, section
5.
One finds that
$$
  \dim H^{2,1}_\omega = 4, \quad
  \dim H^{1,2}_\omega = 1, \qquad
  \dim H^{2,1}_{\bar\omega} = 1, \quad
  \dim H^{1,2}_{\bar\omega} = 4,
$$
and from the Hodge-Riemann bilinear relations one finds that $h$
has signature $(4,1)$.

Now let $\phi$ be a generator of the one-dimensional space
$H^{2,1}_{\bar\omega}$  and let $\gamma_1 \commadots \gamma_5$
be a standard basis of $H^3(T,\Z)$ considered as an $\EE$-module.
By this we mean that the $\gamma_k$ are orthogonal and that
$h(\gamma_1,\gamma_1) = -1$ and $h(\gamma_k,\gamma_k) = 1$ for $k > 1$.
Let $v(\phi,\gamma)$ be the vector in $\C^5$ with components
$$
  v_k = \int_{\gamma_k} \phi .
$$
One verifies that $h(v,v) < 0$ where now $h$ is the hermitian form $-|v_1|^2
+ |v_2|^2 + \cdots + |v_5|^2$.  Thus the line generated by $v(\phi,\gamma)$
defines
a point in $B^4 \subset \P^4$, and one checks that $v(\phi,\gamma) \not\in
\HH$.   By
well-known constructions (the work of Griffiths), the period vector defines  a
holomorphic
map from the universal cover of $M$ to the ball which transforms according to
the
projectivized monodromy representation for marked cyclic cubic threefolds. The
proof that
$\Gamma_0$ is the  projective monodromy group relies on the work of Libgober
\cite{6} and the first author \cite{1}. Thus our construction yields a period
map $f: M \map
(B^4 - \HH)/\Gamma_0$.

\section{3. Properties of the period mapping}

We must now show that $f$ is bijective.  For injectivity, consider once
again the period vector $v(\phi,\gamma)$.  The vectors $\gamma_k$ can be
decomposed into
eigenvectors $\gamma_k'$ and $\gamma_k''$ for $\sigma$, with eigenvalues
$\omega$ and $\bar\omega$, respectively.  Let $\hat\gamma_k'$ and
$\hat\gamma_k''$ denote elements of the corresponding dual basis. Because
$\phi$
is an eigenvector with eigenvalue $\bar\omega$, its integral over $\gamma_k'$
vanishes, so that
$$
  \phi = \sum_k  \hat\gamma_k''\int_{\gamma_k''}\phi
       = \sum_k  \hat\gamma_k''\int_{\gamma_k}\phi .
$$
Thus the components of $v(\phi,\gamma)$ determine $\phi$ as an element of
$H^3_{\bar\omega}$.  Consequently the line $\C v(\phi,\gamma)$ determines
the complex Hodge structure $H^3_{\bar\omega}$.  Viewing the Hodge components
of
$H^3_{\bar\omega}$ as subspaces of $H^3(T,\C)$, we may take their conjugates
to determine the complex Hodge structure $H^3_\omega$.  These two complex Hodge
structures determine the Hodge structure on $H^3(T,\Z)$.  Thus, by the Torelli
theorem of Clemens-Griffiths \cite{4}, the period vector
$v(\phi,\gamma)$ determines the cubic threefold $T$ up to isomorphism.
It remains to show that $T$, which perforce is a cyclic cubic threefold,
determines its ramification locus uniquely.  This follows from the fact that
the locus in question is a planar component of the Hessian surface.

To prove surjectivity we first consider a smooth compactification $C$ of $M$ by
a normal
crossing divisor $D$, e.g., the one given by Naruki \cite{9}, as well as the
Satake
compactification $\overline{B^4/\Gamma_0}$, obtained by adding fourty points,
the ``cusps,''
each corresponding to a null point of $P(\F_3^5, q)$.  By well-known
results  \cite{2} in complex variable theory, the period map
has a holomorphic extension to a map $\bar f$ from $C$ to the Satake
compactification.  Since $C$ is compact,
$\bar f$ is open, and $\overline{B^4/\Gamma_0}$ is connected, we conclude that
$\bar f$ is surjective.

\section{4. Boundary components}

To pass from surjectivity of $\bar f$ to surjectivity of $f$, we must show that
$\bar f$
maps the compactifying divisor $D$ to the complement of $(B^4 - \HH)/\Gamma_0$
in the
Satake compactification. To this end write $D$ as a sum of irreducible
components,
$
  D = \bigcup D_i' \cup \bigcup D''_j
$,
where $D_i'$ parametrizes nodal cubic surfaces via the map to the geometric
invariant theory
compactification of the moduli space of smooth cubics, and where in the same
way the $D_j''$ parametrize cubics with an $A_2$ singularity.

Now consider a one-parameter family of cubic surfaces with smooth total
space acquiring a node.  Its local equation near the node has the form $x^2 +
y^2
+ z^2 = t$ and the corresponding family of cyclic cubic  threefolds has the
form
$x^2 + y^2 + z^2 + w^3 = t$.  The local monodromy of the latter has order six,
its
eigenvalues are primitive sixth roots of unity, and the space of vanishing
cycles is two-dimensional.  (These facts are well-known and the relevant
literature and
arguments are summarized in \cite{3}, section 6).  From
\cite{7} we conclude that coefficients of the period vector on vanishing cycles
are of the
form
$
  A(t) t^{1/6} + B(t)t^{5/6}
$
where $A$ and $B$ are holomorphic.  Now the space of vanishing cycles is
invariant under the action of $\sigma$ and so constitutes a rank one
$\EE$-submodule.  One can choose a generator $\delta$ for it
so that $h(\delta,\delta) = 1$, and then one has
$$
    \lim_{t \to 0} \int_\delta \phi = 0.
$$
Thus the limiting value of the period vector lies in the orthogonal
complement of $\delta$.  In other words, $\bar f(D'_i)$ lies in $\HH/\Gamma_0$,
as required.

Consider next a one-parameter family of cubic surfaces with smooth total space
whose central
fiber acquires an $A_2$ singularity.  Its local equation is $x^2 + y^2 + z^3 =
t$
and the corresponding family of cyclic cubic threefolds has local equation
$x^2 + y^2 + z^3  + w^3 = t$.  In this case the local monodromy is of infinite
order.
After replacing $t$ by $t^3$ one finds an expansion of the form
$ \phi(t) = A(t)(\log t)\,\hat\gamma + \hbox{(terms bounded in $t$})$, where
$A(0) \ne 0$
and where $\hat\gamma$ is an integer cohomology class which is isotropic for
$h$.
Consequently the line $\C\phi(t)$ converges to the isotropic line
$\C\hat\gamma$ as $t$
converges to zero, hence converges to a cusp in the Satake compactification.

\section{6. The corollaries}

Finally, we comment on the collaries.  Part (a) is immediate.  For part
(b) let $K$ denote the kernel of the map $\pi_1(M) \map \Gamma_0$.  Then $K$ is
isomorphic
to the fundamental group of $B^4 - \HH$ and it is easy to see that its
abelianization
is not finitely generated.  We remark that $K$ is not free: there are many
sets of commuting elements corresponding to normal crossings of $\HH$.
For (c) we note first that for lattices in semisimple Lie
group of real rank greater than one, the results of Margulis \cite{8}
imply finite generation of all normal subgroups.  The rank one case can be
treated separately, as was shown to us by Michael Kapovich.

\bibliography

\bi{1} D. Allcock, New complex and quaternion-hyperbolic reflection groups,
submitted,
http://www.math.utah.edu/\~{}allcock

\bi{2} A. Borel, Some metric properties of arithmetic quotients of symmetric
spaces and an
extension theorem. Collection of articles dedicated to S. S. Chern and D. C.
Spencer on
their sixtieth birthdays, J. Differential Geometry {\bf 6 } (1972), 543--560.

\bi{3} J. Carlson and D. Toledo, Discriminant complements and kernels of
monodromy
representations, 22 pp, submitted, http://www.math.utah.edu/\~{}carlson

\bi{4} C. H. Clemens and P. A. Griffiths, the
intermediate Jacobian of the cubic threefold, Ann. of Math. {\bf 95 } (1972),
281--356.

\bi{5}  P.A. Griffiths, On the periods of certain rational
integrals: I and II,  Ann. of Math. {\bf 90} (1969), 460-541.

\bi{6} A. Libgober, On the fundamental group of the space of cubic
surfaces,  Math. Zeit. {\bf 162} (1978), 63--67.

\bi{8} B. Malgrange, Int\'egrales asymptotiques et monodromie,
Ann. Sci. Ecole Norm. Sup., ser. 4, tome { \bf 7} (1974), 405--430.

\bi{8} G. A. Margulis, Quotient groups of discrete subgroups and measure
theory, Funct. Anal. Appl. {\bf 12 } (1978), 295--305.

\bi{9} I. Naruki, Cross ratio variety as a moduli space
of cubic surfaces, Proc. London Math. Soc. (3) {\bf 45 } (1982), 1-30.

\endbibliography

\end